\documentclass[12pt]{amsart}
\usepackage{amssymb,amsmath,amsthm,mathtools,amsfonts,times,hyperref,mathrsfs,multirow}
\usepackage{pgfplots}
\pgfplotsset{compat=1.15}
\usepackage{mathrsfs,mathdots}
\usetikzlibrary{arrows}
\usepackage{enumerate}
\usepackage{graphicx}
\usepackage{array}
\usepackage{ytableau}
\usepackage{pstricks,pst-node,graphicx}
\usepackage{tikz,varwidth}
\usepackage{setspace}
\usepackage{float}
\usetikzlibrary{calc}
\usepackage[super]{nth}
\usepackage{tikz}
\usepackage{extarrows}
\usepackage{pgfplots}
\pgfplotsset{compat=1.15}
\usetikzlibrary{arrows}

\newtheorem{theorem}{Theorem}[section]

\theoremstyle{definition}

\theoremstyle{remark}
\newtheorem{remark}{Remark}
\numberwithin{equation}{section}

\allowdisplaybreaks

\newcommand\numberthis{\addtocounter{equation}{1}\tag{\theequation}}
\newcommand\restr[2]{{
  \left.\kern-\nulldelimiterspace 
  #1 
  \littletaller 
  \right|_{#2} 
  }}
\newcommand{\littletaller}{\mathchoice{\vphantom{\big|}}{}{}{}}  

\begin{document}

\title[]
{On partitions with bounded largest part and fixed integral GBG-rank modulo primes}

\author{Alexander Berkovich}
\address{Department of Mathematics, University of Florida, Gainesville
FL 32611, USA}
\email{alexb@ufl.edu}
\author{Aritram Dhar}
\address{Department of Mathematics, University of Florida, Gainesville
FL 32611, USA}
\email{aritramdhar@ufl.edu}

\dedicatory{Dedicated to George E. Andrews on the occasion of his \nth{85} birthday}

\date{\today}



\begin{abstract}
In 2009, Berkovich and Garvan introduced a new partition statistic called the GBG-rank modulo $t$ which is a generalization of the well-known BG-rank. In this paper, we use the Littlewood decomposition of partitions to study partitions with bounded largest part and fixed integral value of GBG-rank modulo primes. As a consequence, we obtain new elegant generating function formulas for unrestricted partitions, self-conjugate partitions, and partitions whose parts repeat a finite number of times.
\end{abstract}
\maketitle

\section{Introduction and Main Results}\label{s1}
A partition is a non-increasing finite sequence $\pi = (\lambda_1,\lambda_2,\lambda_3,\ldots)$ of positive integers where $\lambda_i$'s are called the parts of $\pi$. We denote the number of parts of $\pi$ by $\#(\pi)$ and the largest part of $\pi$ by $l(\pi)$. The norm of $\pi$, denoted by $|\pi|$, is defined as\\ $$|\pi| = \sum\limits_{i\ge 1}\lambda_i.$$\\ We say that $\pi$ is a partition of $n$ if $|\pi| = n$. We may also write a partition $\pi$ in terms of its frequency of parts as\\ $$\pi = (1^{f_1},2^{f_2},3^{f_3},\ldots),$$\\ where $f_i$ is the frequency of the part $i$. Let $\mathcal{P}$ denote the set of all partitions and $\mathcal{D}$ denote the set of all partitions into distinct parts.\\\par The Young diagram of $\pi$ is a convenient way of representing $\pi$ graphically wherein the parts of $\pi$ are depicted as rows of unit squares which are called cells. Given the Young diagram of $\pi$, we label a cell in the $i$th row and $j$th column by the least non-negative integer $\equiv j-i\ (\textrm{mod}\ t)$. The resulting diagram is called a $t$-residue diagram of $\pi$ \cite{Jam-Kerb81}. We can also label cells in the infinite column $0$ and the infinite row $0$ in the same fashion and call the resulting diagram the extended $t$-residue diagram of $\pi$ \cite{Gar-Kim-Stan90}. For each partition $\pi$ and positive integer $t$, we can associate the $t$-dimensional vector\\ $$\Vec{r} = \Vec{r}(\pi,t) = (r_0(\pi,t), r_1(\pi,t), \ldots, r_{t-1}(\pi,t))$$\\ where\\ $$r_i(\pi,t) = r_i,\quad 0\le i\le t-1$$\\ is the number of cells labeled $i$ in the $t$-residue diagram of $\pi$. For example, Figure \ref{fig1} below depicts the $3$-residue diagram for the partition $\pi = (10,7,4,3)$.
\begin{figure}[H]
\centering
\ytableaushort
{{0}{1}{2}{0}{1}{2}{0}{1}{2}{0},{2}{0}{1}{2}{0}{1}{2},{1}{2}{0}{1},{0}{1}{2}}
* {10,7,4,3}
\par\quad\\\par\quad\\
\caption{$3$-residue diagram of the partition $\pi = (10,7,4,3)$}
\label{fig1}
\end{figure}\par In \cite{Ber-Gar08} and \cite{Ber-Gar09}, Berkovich and Garvan defined a partition statistic of $\pi$\\ 
\begin{equation}\label{eq11}
    \text{GBG-rank}(\pi,t) := \sum\limits_{j=0}^{t-1}r_j(\pi,t)\omega^j_t,
\end{equation}\\ where $\omega_t := e^{\frac{2\pi\iota}{t}}$ is a $t$th root of unity and $\iota = \sqrt{-1}$.\\\par We call the statistic in \eqref{eq11} the GBG-rank of $\pi$ mod $t$ and denote it by $\text{GBG}^{(t)}(\pi)$. For example, in Figure \ref{fig1}, for the partition $\pi = (10, 7, 4, 3)$, $\text{GBG}^{(3)}(\pi) = r_0 + r_1\omega_3 + r_2\omega^2_3 = 8 + 8\omega_3 + 8\omega^2_3 = 0$.\\\par The special case $t = 2$ is called the BG-rank of a partition $\pi$ defined as\\
\begin{equation}\label{eq12}
    \text{BG-rank}(\pi) := r_0(\pi,2) - r_1(\pi,2) = r_0 - r_1.
\end{equation}\\ It has been studied extensively in \cite{Ber-Gar08} and \cite{Ber-Yesil08}.\\\par Let $L,m,n$ be non-negative integers. We now recall some notations from the theory of $q$-series that can be found in \cite{And98}.
\begin{align*}
    (a)_L = (a;q)_L &:= \prod_{k=0}^{L-1}(1-aq^k),\\
    (a)_{\infty} = (a;q)_{\infty} &:= \lim_{L\rightarrow \infty}(a)_L\,\,\text{where}\,\,\lvert q\rvert<1.
\end{align*}\\
\par We define the $q$-binomial (Gaussian) coefficient as\\
\begin{align*}
    \left[\begin{matrix}m\\n\end{matrix}\right]_q := \Bigg\{\begin{array}{lr}
        \dfrac{(q)_m}{(q)_n(q)_{m-n}}\quad\text{for } m\ge n\ge 0,\\
        0\qquad\qquad\quad\text{otherwise}.\end{array}
\end{align*}\\
\begin{remark}\label{rmk1}
For $m, n\ge 0$, $\left[\begin{matrix}m+n\\n\end{matrix}\right]_q$ is the generating function for partitions into at most $n$ parts each of size at most $m$ (see \cite[Chapter $3$]{And98}). Also, note that\\ $$\lim_{m\longrightarrow\infty}\left[\begin{matrix}m+n\\n\end{matrix}\right]_q = \dfrac{1}{(q)_n},$$\\ where $1/(q)_n$ is the generating function for partitions into at most $n$ parts.\\ 
\end{remark}
\par If $\tilde{B}_N(k,q)$ denotes the generating function for the number of partitions into parts less than or equal to $N$ with BG-rank equal to $k$, then Berkovich and Uncu \cite[Theorem $3.2$]{Ber-Unc16} showed that for any non-negative integer $N$ and any integer $k$,\\
\begin{equation}\label{eq13}
    \tilde{B}_{2N+\nu}(k,q) = \dfrac{q^{2k^2-k}}{(q^2;q^2)_{N+k}(q^2;q^2)_{N+\nu-k}},
\end{equation}\\
where $\nu\in\{0,1\}$.\\\par Recently, Dhar and Mukhopadhyay \cite{Dha-Mukh23} asked for a direct combinatorial proof of \eqref{eq13}. In this paper, we answer their question in the affirmative and generalize \eqref{eq13} as Theorem \ref{thm11}. For instance, for any odd prime $t$, the generating function for partitions $\pi$ having largest part $l(\pi)\le tN$ and $\text{GBG}^{(t)}(\pi) = k\in\mathbb{Z}$ is\\
\begin{align*}
\dfrac{q^{tk^2-(t-1)k}}{(q^t;q^t)_{N-k}(q^t;q^t)_N^{t-2}(q^t;q^t)_{N+k}},    
\end{align*}\\ 
which is the $\nu = 0$ special case of the generating function in Theorem \ref{thm11} below.\\\par Let $G_{N,t}(k,q)$ denote the generating function for the number of partitions into parts less than or equal to $N$ with GBG-rank mod $t$ equal to $k$, i.e.,\\
$$G_{N,t}(k,q) := {\displaystyle\sum\limits_{\substack{\pi\in\mathcal{P} \\ \text{GBG}^{(t)}(\pi) = k \\ l(\pi)\le N}}q^{|\pi|}}.$$\\ We then have the following new result.\\
\begin{theorem}\label{thm11}
For any prime $t$, a non-negative integer $N$, and any integer $k$, we have
\begin{equation}\label{eq14}
G_{tN+\nu,t}(k,q) = \dfrac{q^{tk^2-(t-1)k}}{\displaystyle\prod_{i=0}^{t-1}(\tilde{q};\tilde{q})_{N+\lceil\frac{\nu-i}{t}\rceil-k\delta_{i,0}+k\delta_{i,t-1}}},    
\end{equation}
where $\tilde{q} := q^t$, and $\nu\in\{0,1,2,\ldots,t-1\}$.\\ 
\end{theorem}
Above and throughout the paper, $\lfloor x\rfloor$ denotes the floor function, i.e, the greatest integer less than or equal to $x$, $\lceil x\rceil$ denotes the ceil function, i.e., the least integer greater than or equal to $x$, and $\delta_{i,j}$ denotes the Kronecker delta function, i.e., $\delta_{i,j}$ is equal to $1$ if $i = j$ otherwise $0$.\\\par Recall \cite{And98} that the \textit{conjugate} of a partition $\pi$, denoted by $\pi^{\prime}$, is associated with the Young diagram obtained by reflecting the diagram of $\pi$ across the main diagonal. We say that $\pi$ is \textit{self-conjugate} if $\pi = \pi^{\prime}$. Let $\mathcal{SCP}$ denote the set of all self-conjugate partitions.\\
\begin{remark}\label{rmk2}
    $(-q;q^2)_L$ is the generating function for partitions into district odd parts having the largest part at most $2L - 1$. Since self-conjugate partitions are in bijection with distinct odd part partitions, $(-q;q^2)_L$ is also the generating function for self-conjugate partitions with the number of parts at most $L$ (see \cite{And98}).\\
\end{remark}
Let $GSC_{N,t}(k,q)$ denote the generating function for the number of self-conjugate partitions into parts less than or equal to $N$ with GBG-rank mod $t$ equal to $k$, i.e.,\\ 
$$GSC_{N,t}(k,q) := {\displaystyle\sum\limits_{\substack{\pi\in\mathcal{SCP} \\ \text{GBG}^{(t)}(\pi) = k \\ l(\pi)\le N}}q^{|\pi|}}.$$\\ We then have the following new results.\\
\begin{theorem}\label{thm12}
For any non-negative integer $N$ and any integer $k$, we have\\
\begin{equation}\label{eq15}
GSC_{2N+\nu,2}(k,q) = q^{2k^2-k}\left[\begin{matrix}2N+\nu\\N+k\end{matrix}\right]_{q^4}    ,
\end{equation}
where $\nu\in\{0,1\}.$\\
\end{theorem}
\begin{theorem}\label{thm13}
For any odd prime $t$, a non-negative integer $N$, and any integer $k$, we have\\
\begin{align*}
GSC_{tN+\nu,t}(k,q) &= q^{tk^2-(t-1)k}(-q^t;q^{2t})_{N+\big\lceil\frac{\nu-\frac{t-1}{2}}{t}\big\rceil}\left[\begin{matrix}2N+\lceil\frac{\nu}{t}\rceil\\N+k\end{matrix}\right]_{q^{2t}}\\&\qquad\times{\displaystyle\prod_{i=1}^{\frac{t-3}{2}}\left[\begin{matrix}2N+\lceil\frac{\nu-i}{t}\rceil+\lfloor\frac{\nu+i}{t}\rfloor\\N+\lfloor\frac{\nu+i}{t}\rfloor\end{matrix}\right]_{q^{2t}}},\numberthis\label{eq16}\end{align*}
where $\nu\in\{0,1,2,\ldots,t-1\}$.\\
\end{theorem}
\par If $B_N(k,q)$ denotes the generating function for the number of partitions into distinct parts less than or equal to $N$ with BG-rank equal to $k$, then Berkovich and Uncu \cite[Theorem $3.1$]{Ber-Unc16} showed that for any non-negative integer $N$ and any integer $k$,\\
\begin{equation}\label{eq17}
    B_{2N+\nu}(k,q) = q^{2k^2-k}\left[\begin{matrix}2N+\nu\\N+k\end{matrix}\right]_{q^2},
\end{equation}\\
where $\nu\in\{0,1\}$.\\\par Now, we generalize \eqref{eq17} as below.\\\par let $\tilde{G}_{N,t}(k,q)$ denote the generating function for the number of partitions into parts repeating no more than $t-1$ times and less than or equal to $N$ with GBG-rank mod $t$ equal to $k$, i.e.,\\
$$\tilde{G}_{N,t}(k,q) := {\displaystyle\sum\limits_{\substack{\pi = (1^{f_1},2^{f_2},\ldots,N^{f_N})\in\mathcal{P} \\ f_i\le t-1, 1\le i\le N \\ \text{GBG}^{(t)}(\pi) = k}}q^{|\pi|}}.$$\\ We then have the following new result.\\
\begin{theorem}\label{thm14}
For any prime $t$, a non-negative integer $N$, and any integer $k$, we have
\begin{equation}\label{eq18}
\tilde{G}_{tN+\nu,t}(k,q) = \dfrac{q^{tk^2-(t-1)k}(\tilde{q};\tilde{q})_{tN+\nu}}{\displaystyle\prod_{i=0}^{t-1}(\tilde{q};\tilde{q})_{N+\lceil\frac{\nu-i}{t}\rceil-k\delta_{i,0}+k\delta_{i,t-1}}},    
\end{equation}
where $\tilde{q} := q^t$, and $\nu\in\{0,1,2,\ldots,t-1\}$.\\
\end{theorem}
The rest of the paper is organized as follows. In Section \ref{s2}, we introduce preliminary background on rim hooks, $t$-cores, Littlewood decomposition, and $t$-quotients of partitions. In Section \ref{s3}, we prove Theorems \ref{thm11} - \ref{thm14} using Littlewood decomposition of partitions. In Section \ref{s4}, we conclude with some interesting observations.\\

\section{Background on Littlewood decomposition of partitions}\label{s2}
We now recall the notions of rim hook and $t$-core \cite{Jam-Kerb81}. If a cell of $\pi$ shares a vertex or edge with the rim of the diagram of $\pi$, we call this cell a rim cell of $\pi$. A connected collection of rim cells of $\pi$ is called a rim hook if (Young diagram of $\pi$)$\setminus$(rim hook) represents a legitimate partition. We then call a partition a $t$-core, denoted $\pi_{\text{$t$-core}}$, if its Young diagram has no rim hooks of length $t$ \cite{Jam-Kerb81}. Throughout the paper we will denote a $t$-core by $\pi_{\text{$t$-core}}$. Any partition $\pi$ has a uniquely determined $t$-core which we also denote by $\pi_{\text{$t$-core}}$. This partition $\pi_{\text{$t$-core}}$ is called the $t$-core of $\pi$. One can obtain $\pi_{\text{$t$-core}}$ from $\pi$ by the successive removal of rim hooks of length $t$. The $t$-core $\pi_{\text{$t$-core}}$ is independent of the way in which the hooks are removed.\\\par We now recall some definitions from \cite{Gar-Kim-Stan90}. A region $r$ in the extended $t$-residue diagram of $\pi$ is the set of all cells $(i,j)$ satisfying $t(r-1)\le j-i< tr$. A cell of $\pi$ is called exposed if it is at the end of a row in the extended $t$-residue diagram of $\pi$. One can construct $t$ infinite binary words $W_0,W_1,\ldots,W_{t-1}$ of two letters $N, E$ as follows: The $r$th letter of $W_i$ is $E$ if there is an exposed cell labeled $i$ in the region $r$, otherwise the $r$th letter of $W_i$ is $N$. It is then easy to see that the word set $\{W_0,W_1,\ldots,W_{t-1}\}$ fixes $\pi$ uniquely. For example, the three bi-infinite words $W_0,W_1,W_2$ for the partition $(10,7,4,3)$ in Figure \ref{fig1} are as follows:\\\par
Region : $\cdots\cdots$\quad$-3$\qquad$-2$\qquad$-1$\qquad$0$\qquad$1$\qquad$2$\qquad$3$\qquad$4$\qquad$5$\quad$\cdots\cdots$\\\par
\quad\,\,\,\,$W_0$ : $\cdots\cdots$\quad\,\,$E$\qquad\,\,\,$E$\qquad\,\,\,$E$\quad\,\,\,\,$N$\quad\,\,\,\,$N$\quad\,\,\,\,$N$\quad\,\,\,\,$N$\quad\,\,\,\,$E$\quad\,\,\,\,$N$\,\,\,\,$\cdots\cdots$\\\par
\quad\,\,\,\,$W_1$ : $\cdots\cdots$\quad\,\,$E$\qquad\,\,\,$E$\qquad\,\,\,$E$\quad\,\,\,\,$N$\quad\,\,\,\,$E$\quad\,\,\,\,$N$\quad\,\,\,\,$N$\quad\,\,\,\,$N$\quad\,\,\,\,$N$\,\,\,\,$\cdots\cdots$\\\par
\quad\,\,\,\,$W_2$ : $\cdots\cdots$\quad\,\,$E$\qquad\,\,\,$E$\qquad\,\,\,$N$\quad\,\,\,\,$E$\quad\,\,\,\,$N$\quad\,\,\,\,$E$\quad\,\,\,\,$N$\quad\,\,\,\,$N$\quad\,\,\,\,$N$\,\,\,\,$\cdots\cdots$\\\par Let $\mathcal{P}_\text{$t$-core}$ denote the set of all $t$-cores. There is a well-known bijection\\ $$\phi_1 : \mathcal{P}\longrightarrow\mathcal{P}_{\text{$t$-core}}\times\mathcal{P}\times\mathcal{P}\times\mathcal{P}\times\cdots\times\mathcal{P}$$\\ due to Littlewood \cite{Little51}\\ $$\phi_1(\pi) = (\pi_{\text{$t$-core}},(\Hat{\pi}_0,\Hat{\pi}_1,\ldots,\Hat{\pi}_{t-1}))$$\\ such that\\ $$|\pi| = |\pi_{\text{$t$-core}}| + t\sum\limits_{i=0}^{t-1}|\Hat{\pi}_i|.$$\\\par The vector partition $(\Hat{\pi}_0,\Hat{\pi}_1,\ldots,\Hat{\pi}_{t-1})$ is called the $t$-quotient of $\pi$ and is denoted by $\pi_{\text{$t$-quotient}}$. Let $\mathcal{P}_{\text{$t$-quotient}}$ be the set of all $t$-quotients.\\\par We will now define the vector\\ $$\Vec{n} = \Vec{n}(\pi,t) = (n_0(\pi,t),n_1(\pi,t),\ldots,n_{t-1}(\pi,t)) = (n_0,n_1,\ldots,n_{t-1})$$\\ where for $0\le i\le t-2$\\ $$n_i = r_i - r_{i+1},$$\\ and\\ $$n_{t-1} = r_{t-1} - r_0.$$\\\par Note that\\ $$\sum\limits_{i=0}^{t-1}n_i = 0.$$\\\par Moreover, it can be shown \cite{Gar-Kim-Stan90} that\\ $$|\pi_{\text{$t$-core}}| = \frac{t}{2}\Vec{n}\cdot\Vec{n} + \Vec{b}_t\cdot\Vec{n},$$\\ where $\Vec{b}_t := (0,1,2,\ldots,t-1)$.\\\par For $0\le i\le t-1$, define $\chi_i(\pi,t)$ to be the largest region in the extended $t$-residue diagram of $\pi$ where the cell labeled $i$ is exposed. In \cite{Ber-Gar08}, Berkovich and Garvan observed that\\
\begin{align*}
    \chi_i(\pi,t) = v_i + n_i(\pi,t),\numberthis\label{eq21}\\
\end{align*} where $v_i$ is the number of parts in the $i$th component of the $t$-quotient of $\pi$ and $n_i(\pi,t)$ is the $i$th component of $\Vec{n}(\pi,t)$. For example, for the partition $\pi = (10,7,4,3)$ in Figure \ref{fig1}, from its three bi-infinite words $W_0$, $W_1$, and $W_2$ as shown above, we have $\chi_0(\pi,3) = 4$, $\chi_1(\pi,3) = 1$, and $\chi_2(\pi,3) = 2$. Note that $\Vec{n}(\pi,t) = (0,0,0)$. Therefore, using \eqref{eq21}, we observe that the $0$-component of $\pi_{\text{$3$-quotient}}$ has $4$ parts, the $1$-component of $\pi_{\text{$3$-quotient}}$ has $1$ part, and the $2$-component of $\pi_{\text{$3$-quotient}}$ has $2$ parts.\\\par Note that if $\pi^{\prime}$ is the conjugate of $\pi\in\mathcal{P}$ whose Littlewood decomposition is\\ $$\phi_1(\pi) = (\pi_{\text{$t$-core}},\pi_{\text{$t$-quotient}})$$\\ where\\ $$\pi_{\text{$t$-quotient}} = (\Hat{\pi}_0,\Hat{\pi}_1,\ldots,\Hat{\pi}_{t-1}),$$\\ then under conjugation, we have\\ $$\phi_1(\pi^{\prime}) = (\pi^{\prime}_{\text{$t$-core}},\pi^{\prime}_{\text{$t$-quotient}})$$\\ where $\pi^{\prime}_{\text{$t$-core}}$ is the conjugate of $\pi_{\text{$t$-core}}$ and\\ $$\pi^{\prime}_{\text{$t$-quotient}} = (\Hat{\pi}^{\prime}_{t-1},\Hat{\pi}^{\prime}_{t-2},\ldots,\Hat{\pi}^{\prime}_0)$$\\ is the conjugate of $\pi_{\text{$t$-quotient}} = (\Hat{\pi}_0,\Hat{\pi}_1,\ldots,\Hat{\pi}_{t-1})$. Please examine \cite[p.19, Proposition 3.12]{Bru-Nath19} for a brief explanation of the above relation. Also, observe that under conjugation,\\ $$\Vec{n}(\pi,t) = (n_0,n_1,\ldots,n_{t-2},n_{t-1})$$\\ becomes\\ $$\Vec{\tilde{n}} = \Vec{n}(\pi^{\prime},t) = (-n_{t-1},-n_{t-2},\ldots,-n_1,-n_0).$$\\\par Please refer to \cite{Gar-Kim-Stan90} for the above observation. \\\par Now, suppose $\pi\in\mathcal{SCP}$. Then, we have\\ $$\phi_1(\pi) = (\pi_{\text{$t$-core}},\pi_{\text{$t$-quotient}})$$\\ where\\ $$\pi_{\text{$t$-core}} = \pi^{\prime}_{\text{$t$-core}}$$\\ and\\ $$\pi_{\text{$t$-quotient}} = (\Hat{\pi}_0,\Hat{\pi}_1,\Hat{\pi}_2,\ldots,\Hat{\pi}_{t-1})$$\\ with\\ $$\Hat{\pi}_i = \Hat{\pi}^{\prime}_{t-1-i}\,,\,\, 0\le i\le t-1.$$\\\par Note that for $t$ being odd,\\ $$\Hat{\pi}_{\frac{t-1}{2}} = \Hat{\pi}^{\prime}_{\frac{t-1}{2}}.$$\\

\section{Proofs of Theorems \ref{thm11}-\ref{thm14}}\label{s3}
In this section, we provide new combinatorial proofs of our main results in Section \ref{s1}.
\subsection{Proof of Theorem \ref{thm11}}\label{ss31}
Let $\pi = (\lambda_1,\lambda_2,\lambda_3,\ldots)\in\mathcal{P}$ be a partition with $\lambda_1\le tN+\nu$ where $0\le\nu\le t-1$ and consider the Littlewood decomposition of $\pi$\\ $$\phi_1(\pi) = (\pi_{\text{$t$-core}},\pi_{\text{$t$-quotient}})$$\\ where $\pi_{\text{$t$-quotient}} = (\Hat{\pi}_0,\Hat{\pi}_1,\ldots,\Hat{\pi}_{t-1})$.\\\par Now, let $t$ be any prime and $\text{GBG}^{(t)}(\pi) = k$ be an integer. One can then easily verify that $\text{GBG}^{(t)}(\pi) = k$ is an integer if and only if $$r_0(\pi,t) = k + r$$ and $$r_i(\pi,t) = r$$ for any $r\in\mathbb{Z}^{+}\cup\{0\}$ and $1\le i\le t-1$. We leave it as an exercise to the reader.\\\par Therefore, we have\\ $$\Vec{n}(\pi,t) = (k,0,0,\ldots,0,0,-k),$$\\ where $n_0(\pi,t) = k$, $n_{t-1}(\pi,t) = -k$, and $n_i(\pi,t) = 0$ for $1\le i\le t-2$.\\\par Thus, we have\\
\begin{align*}
    |\pi_{\text{$t$-core}}| &= \frac{t}{2}\cdot 2k^2 + (t-1)\cdot (-k)\\
    &= tk^2 - (t-1)k.\numberthis\label{eq31}\\
\end{align*}\par Now, following Figure \ref{fig2} below, we consider two cases regarding $\nu$ which are as follows:\\
\vspace{-1cm}
\begin{figure}[H]
\definecolor{ccqqqq}{rgb}{1.0,0,0}
\begin{tikzpicture}[line cap=round,line join=round,>=triangle 45,x=0.7cm,y=0.7cm]
\clip(-3.0,-15.0) rectangle (28.0,3.0);
\draw [line width=1pt] (0,0)-- (0,-7);
\draw [line width=1pt] (0,0)-- (-1,1);
\draw [line width=1pt] (3,0)-- (2,1);
\draw [line width=1pt] (0,-3)-- (-1,-2);
\draw [line width=1pt] (0,0)-- (3,-3);
\draw [line width=1pt] (0,-3)-- (2,-5);
\draw [line width=1pt] (3,0)-- (5,-2);
\draw [line width=1pt] (0,0)-- (6,0);
\draw (6.0,0.2) node[anchor=north west] {$\ldots\ldots$};
\draw [line width=1pt] (8,0)-- (7.009613178261587,1.0042800934430731);
\draw (8.2,0.3) node[anchor=north west] {\tiny{$0$}};
\draw (8.6,0.3) node[anchor=north west] {\tiny{$1$}};
\draw (9.0,0.3) node[anchor=north west] {\tiny{$2$}};
\draw (9.4,0.3) node[anchor=north west] {\tiny{$3$}};
\draw [line width=1pt] (8,0)-- (10.726343675686474,-2.733810924081948);
\draw (10.3,-2.0) node[anchor=north west] {\tiny{$0$}};
\draw (10.3,-1.6) node[anchor=north west] {\tiny{$1$}};
\draw (10.3,-1.2) node[anchor=north west] {\tiny{$2$}};
\draw (10.7,-1.2) node[anchor=north west] {\tiny{$3$}};
\draw (11.1,-1.2) node[anchor=north west] {\tiny{$4$}};
\draw (11.1,-0.8) node[anchor=north west] {\tiny{$5$}};
\draw (12.2,0.3) node[anchor=north west] {\tiny{$\nu-1$}};
\draw (12.6,0.9) node[anchor=north west] {\tiny{$\nu$}};
\draw (13.3,1.0) node[anchor=north west] {\tiny{$\nu+1$}};
\draw [line width=1pt] (13.61588642187619,0)-- (17,0);
\draw [line width=1pt] (17,0)-- (19,-2);
\draw [line width=1pt] (17,0)-- (19,0);
\draw (16.4,1.0) node[anchor=north west] {\tiny{$t-1$}};
\draw (9.8,0.2) node[anchor=north west] {$\ldots\ldots\ldots$};
\draw (14.3,0.9) node[anchor=north west] {$\ldots\ldots.\,.$};
\draw (11.4,0.2) node[anchor=north west] {$\iddots$};
\draw (-0.6,1.0) node[anchor=north west] {\tiny{$0$}};
\draw (-0.2,1.0) node[anchor=north west] {\tiny{$1$}};
\draw (0.2,0.9) node[anchor=north west] {$\ldots$};
\draw (1.1,1.0) node[anchor=north west] {\tiny{$t-1$}};
\draw [color=ccqqqq](0.9,-1.8) node[anchor=north west] {\small{$0$}};
\draw [color=ccqqqq](2.6,-0.7) node[anchor=north west] {\small{$1$}};
\draw [color=ccqqqq](0.3,-4.5) node[anchor=north west] {\small{$-1$}};
\draw [color=ccqqqq](4.7,-0.2) node[anchor=north west] {\small{$2$}};
\draw [color=ccqqqq](8.3,-1.1) node[anchor=north west] {\small{$N$}};
\draw [color=ccqqqq](15.0,-0.8) node[anchor=north west] {\small{$N+1$}};
\draw [color=ccqqqq](17.7,-0.3) node[anchor=north west] {\small{$N+2$}};
\end{tikzpicture}
\vspace{-5cm}
\caption{Generic extended $t$-residue diagram of a partition $\pi = (\lambda_1,\lambda_2,\lambda_3,\ldots)\in\mathcal{P}$ with $\lambda_1 = tN+\nu$ where $0\le\nu\le t-1$ (the regions are labeled in red).}
\label{fig2}
\end{figure}
\par\quad\\
\begin{itemize}
    \item \fbox{Case I: $\nu > 0$}\\\par Observe that it is easy to check, from Figure \ref{fig2}, that the cell labeled $i$ may be exposed in the region marked $N+1$ for $0\le i\le \nu-1$ and is not exposed in the region marked $N+1$ for $\nu\le i\le t-1$. Thus, it follows that\\ $$N + 1\ge\chi_i(\pi,t)\,,\quad 0\le i\le \nu-1,$$ and $$N\ge\chi_i(\pi,t)\,,\quad \nu\le i\le t-1.$$\\ Equation \eqref{eq21} then implies that\\
    \begin{equation}\label{eq32}
        N + 1\ge v_i + k\delta_{i,0}\,,\quad 0\le i\le \nu-1,
    \end{equation} and
    \begin{equation}\label{eq33}
        N\ge v_i - k\delta_{i,t-1}\,,\quad \nu\le i\le t-1.
    \end{equation}\\\par
    \item  \fbox{Case II: $\nu = 0$}\\\par Again, from Figure \ref{fig2}, the cell labeled $i$ is not exposed in the region marked $N+1$ for $0\le i\le t-1$. Thus, it follows that\\ $$N\ge\chi_i(\pi,t)\,,\quad 0\le i\le t-1.$$\\ Equation \eqref{eq21} then implies that\\
    \begin{equation}\label{eq34}
        N\ge v_i + k\delta_{i,0} - k\delta_{i,t-1}\,,\quad 0\le i\le t-1.
    \end{equation}\\
\end{itemize}\par 
Hence, combining \eqref{eq32}, \eqref{eq33}, and \eqref{eq34}, we have\\
\begin{equation}\label{eq35}
v_i\le N + \Big\lceil\frac{\nu-i}{t}\Big\rceil - k\delta_{i,0} + k\delta_{i,t-1}\,,\quad 0\le i\le t-1,\\   
\end{equation}\\ where $v_i$ is the number of parts in the $i$th component of the $t$-quotient of $\pi$.\\\par Hence, the required generating function follows from Remark \ref{rmk1}, \eqref{eq31}, and \eqref{eq35} which completes the proof of Theorem \ref{thm11}.\qed\\\par
\begin{remark}\label{rmk3}
Observe that \eqref{eq13} is a special case of Theorem \ref{thm11} where $t = 2$ and so, for $t = 2$, the proof of Theorem \ref{thm11} provides a direct combinatorial proof of \eqref{eq13} as asked by Dhar and Mukhopadhyay in \cite{Dha-Mukh23}.\\
\end{remark}

\subsection{Proof of Theorem \ref{thm12}}\label{ss32}
Suppose $\pi\in\mathcal{SCP}$ be a self-conjugate partition having $l(\pi)\le 2N+1$ and consider the Littlewood decomposition of $\pi$\\ $$\phi_1(\pi) = (\pi_{\text{$2$-core}},\pi_{\text{$2$-quotient}})$$\\ where $\pi_{\text{$2$-quotient}} = (\Hat{\pi}_0,\Hat{\pi}_1)$ and $\Hat{\pi}_1 = \Hat{\pi}^{\prime}_0$.\\\par Now, let $\text{BG-rank}(\pi) = k\in\mathbb{Z}$. It is again easy to verify that $\text{BG-rank}(\pi) = k$ if and only if $$r_0(\pi) = k + r$$ and $$r_1(\pi) = r$$ for any $r\in\mathbb{Z}^{+}\cup\{0\}$ which implies that\\
\begin{align*}
    \Vec{n}(\pi,2) = (k,-k).\\
\end{align*}\par Thus, we have\\
\begin{align*}
    |\pi_{\text{$2$-core}}| &= \frac{2}{2}\cdot 2k^2 + (1)\cdot (-k)\\
    &= 2k^2 - k.\numberthis\label{eq36}
\end{align*}\\\par Now, observe that the cell labeled $0$ may be exposed in the region marked $N+1$ and the cell labeled $1$ is not exposed in the region marked $N+1$. Thus, again using \eqref{eq21}, the following equations hold:
\begin{equation}\label{eq37}
N + 1\ge k + v_0,
\end{equation} which implies that
\begin{equation}\label{eq38}
v_0\le N + 1 - k,
\end{equation} and
\begin{equation}\label{eq39}
N\ge -k + v_1,
\end{equation} which implies that
\begin{equation}\label{eq310}
v_1\le N + k.
\end{equation}\\\par Hence, from \eqref{eq38} and \eqref{eq310}, we have $l(\Hat{\pi}_0)\le N + k$  and $\#(\Hat{\pi}_0)\le N + 1 - k$ which when combined with Remark \ref{rmk1} and \eqref{eq36} gives us the required generating function in Theorem \ref{thm12} for $\nu = 1$. Analogously, one can prove the required generating function for $\nu = 0$ which gives us the complete proof of Theorem \ref{thm12}.\qed\\

\subsection{Proof of Theorem \ref{thm13}}\label{ss33}
Let $\pi = (\lambda_1,\lambda_2,\lambda_3,\ldots)\in\mathcal{SCP}$ be a partition with $\lambda_1\le tN+\nu$ where $0\le\nu\le t-1$ and consider the Littlewood decomposition of $\pi$\\ $$\phi_1(\pi) = (\pi_{\text{$t$-core}},\pi_{\text{$t$-quotient}})$$\\ where $\pi_{\text{$t$-quotient}} = (\Hat{\pi}_0,\Hat{\pi}_1,\ldots,\Hat{\pi}_{t-1})$.\\\par Now, let $t$ be any odd prime and $\text{GBG}^{(t)}(\pi) = k$ be an integer. It is then easy to verify that $\text{GBG}^{(t)}(\pi) = k$ is an integer if and only if $$r_0(\pi,t) = k + r$$ and $$r_i(\pi,t) = r$$ for any $r\in\mathbb{Z}^{+}\cup\{0\}$ and $1\le i\le t-1$.\\\par Therefore, we have\\
\begin{align*}
    \Vec{n}(\pi,t) = (k,0,0,\ldots,0,0,-k),\\
\end{align*} where $n_0(\pi,t) = k$, $n_{t-1}(\pi,t) = -k$, and $n_i(\pi,t) = 0$ for $1\le i\le t-2$.\\\par Thus, we have\\
\begin{align*}
    |\pi_{\text{$t$-core}}| &= \frac{t}{2}\cdot 2k^2 + (t-1)\cdot (-k)\\
    &= tk^2 - (t-1)k.\numberthis\label{eq311}\\
\end{align*}\par
Now, let us consider the pair of partitions $(\Hat{\pi}_0,\Hat{\pi}_{t-1})$ from $\pi_{\text{$t$-quotient}}$ and observe that $\Hat{\pi}_0 = \Hat{\pi}^{\prime}_{t-1}$. Therefore, this pair contributes to the generating function as\\
\begin{align*}
    q^{t(|\Hat{\pi}_0|+|\Hat{\pi}_{t-1}|)} = q^{2t|\Hat{\pi}_0|}.\numberthis\label{eq312}\\
\end{align*}\par Now, by \eqref{eq35}, we have $\#(\Hat{\pi}_0)\le N - k + \lceil\frac{\nu}{t}\rceil$ and $\#(\Hat{\pi}_{t-1})\le N + k$ which is the same as $l(\Hat{\pi}_0)\le N + k$. So, using Remark \ref{rmk1} and \eqref{eq312}, we get the term\\
\begin{align*}
    \left[\begin{matrix}2N+\lceil\frac{\nu}{t}\rceil\\N+k\end{matrix}\right]_{q^{2t}}\numberthis\label{eq313}\\
\end{align*} in the required generating function.\\\par Next, for $1\le i\le (t-3)/2$, we consider the pair of partitions $(\Hat{\pi}_i,\Hat{\pi}_{t-1-i})$ from $\pi_{\text{$t$-quotient}}$ and observe that $\Hat{\pi}_i = \Hat{\pi}^{\prime}_{t-1-i}$. Therefore, this pair contributes to the generating function as\\
\begin{align*}
    q^{t(|\Hat{\pi}_i|+|\Hat{\pi}_{t-1-i}|)} = q^{2t|\Hat{\pi}_i|}.\numberthis\label{eq314}\\
\end{align*}\par Now, by \eqref{eq35}, for $1\le i\le (t-3)/2$, we have $\#(\Hat{\pi}_i)\le N + \lceil\frac{\nu-i}{t}\rceil$ and $\#(\Hat{\pi}_{t-1-i})\le N + \lceil\frac{\nu-(t-1-i)}{t}\rceil = N + \lfloor\frac{\nu+i}{t}\rfloor$ which is the same as $l(\Hat{\pi}_i)\le N + \lfloor\frac{\nu+i}{t}\rfloor$. So, using Remark \ref{rmk1} and \eqref{eq314} followed by taking the product over all values of $i$, we get the term\\
\begin{align*}
    {\displaystyle\prod_{i=1}^{\frac{t-3}{2}}\left[\begin{matrix}2N+\lceil\frac{\nu-i}{t}\rceil+\lfloor\frac{\nu+i}{t}\rfloor\\N+\lfloor\frac{\nu+i}{t}\rfloor\end{matrix}\right]_{q^{2t}}}\numberthis\label{eq315}\\
\end{align*} in the required generating function.\\\par Finally, observe that $\Hat{\pi}_{\frac{t-1}{2}} = \Hat{\pi}^{\prime}_{\frac{t-1}{2}}$ and similarly, by \eqref{eq35}, we have $\#(\Hat{\pi}_{\frac{t-1}{2}})\le N+\big\lceil\frac{\nu-\frac{t-1}{2}}{t}\big\rceil$. Thus, using Remark \ref{rmk2}, we get the term\\
\begin{align*}
    (-q^t;q^{2t})_{N+\big\lceil\frac{\nu-\frac{t-1}{2}}{t}\big\rceil}\numberthis\label{eq316}\\
\end{align*} in the required generating function.\\\par Hence, combining \eqref{eq311}, \eqref{eq313}, \eqref{eq315}, and \eqref{eq316}, we have the desired generating function which completes the proof of Theorem \ref{thm13}.\qed\\

\subsection{Proof of Theorem \ref{thm14}}\label{ss34}
The proof of Theorem \ref{thm14} comes from the combinatorial bijection of extracting parts which repeat $t$ times in succession from the $t$-residue diagram of a partition to form a new pair of partitions.\\\par Let\\ $$\pi = (1^{tq_1+f_1},2^{tq_2+f_2},3^{tq_3+f_3},\ldots,N^{tq_N+f_N}),$$\\ where $0\le f_j\le t-1$ and $q_j\ge 0$ for $1\le j\le N$.\\\par Then, it can be seen that $\pi$ is in one-to-one correspondence with the pair of partitions $(\pi_1,\pi_2)$ where\\ $$\pi_1 = (1^{tq_1},2^{tq_2},3^{tq_3},\ldots,N^{tq_N})$$ and $$\pi_2 = (1^{f_1},2^{f_2},3^{f_3},\ldots,N^{f_N}).$$\par Now, note that\\
\begin{align*}
    \sum\limits_{\pi_1}q^{|\pi_1|} = \dfrac{1}{(q^t;q^t)_N},\numberthis\label{eq317}\\
\end{align*}
and so, we have\\
\begin{align*}
    \sum\limits_{\pi}q^{|\pi|} &= \sum\limits_{\pi_1}q^{|\pi_1|}\cdot\sum\limits_{\pi_2}q^{|\pi_2|}\\
    &= \dfrac{1}{(q^t;q^t)_N}\sum\limits_{\pi_2}q^{|\pi_2|},\numberthis\label{eq318}\\
\end{align*}
where \eqref{eq318} follows from \eqref{eq317}.\\\par Hence, from \eqref{eq318}, we have\\
\begin{align*}
    \sum\limits_{\pi_2}q^{|\pi_2|} = (q^t;q^t)_N\sum\limits_{\pi}q^{|\pi|},\numberthis\label{eq319}\\
\end{align*}
where $\pi_2$ is a partition whose parts repeat no more than $t-1$ times.\\\par Also, observe that $\text{GBG}^{(t)}(\pi) = \text{GBG}^{(t)}(\pi_2) = k$ since removal of parts which repeat $t$ times in succession keeps the GBG-rank mod $t$ value of $\pi$ invariant.\\\par Therefore, using the above observation and replacing $N$ by $tN+\nu$ in \eqref{eq319} for $0\le\nu\le t-1$ along with Theorem \ref{thm11} gives us the desired generating function $\tilde{G}_{tN+\nu,t}(k,q)$.\qed\\\par
\begin{remark}\label{rmk4}
Note that, for $t=2$, we get \cite[eq. (3.2)]{Ber-Unc16} for distinct part partitions, i.e.,\\
$${\displaystyle\sum\limits_{\substack{\pi\in\mathcal{D} \\ \text{BG-rank}(\pi) = k \\ l(\pi)\le 2N+\nu}}q^{|\pi|}} = q^{2k^2-k}\left[\begin{matrix}2N+\nu\\N+k\end{matrix}\right]_{q^2}.$$\\
\end{remark}

\section{Further Observations}\label{s4}
In this section, we will see that Theorem \ref{thm11} can be easily extended to the case where partitions have GBG-rank mod $t$ value equal to $k\omega^j_t$ for any odd prime $t$, any integer $k$ and $1\le j\le t-1$.\\\par Let $\pi = (\lambda_1,\lambda_2,\lambda_3,\ldots)\in\mathcal{P}$ be a partition with $\lambda_1\le tN+\nu$ where $0\le\nu\le t-1$ and consider the Littlewood decomposition of $\pi$\\ $$\phi_1(\pi) = (\pi_{\text{$t$-core}},\pi_{\text{$t$-quotient}})$$\\ where $\pi_{\text{$t$-quotient}} = (\Hat{\pi}_0,\Hat{\pi}_1,\ldots,\Hat{\pi}_{t-1})$.\\\par Now, let $t$ be any odd prime and $\text{GBG}^{(t)}(\pi) = k\omega^j_t$ where $k$ is an integer and $1\le j\le t-1$. One can then easily verify that $\text{GBG}^{(t)}(\pi) = k\omega^j_t$ if and only if\\ $$r_i(\pi,t) = r + k\delta_{i,j}$$\\ for any $r\in\mathbb{Z}^{+}\cup\{0\}$, $0\le i\le t-1$, and $1\le j\le t-1$.\\\par Therefore, we have\\ $$\Vec{n}(\pi,t) = (n_0,n_1,\ldots,n_{t-1}),$$\\ where $n_i(\pi,t) = -k\delta_{i,j-1} + k\delta_{i,j}$ for $0\le i\le t-1$ and $1\le j\le t-1$.\\\par Thus, we have\\
\begin{align*}
    |\pi_{\text{$t$-core}}| = tk^2 + k.\numberthis\label{eq41}
\end{align*}\\\par We then have the following new result.\\
\begin{theorem}\label{thm41}
For any odd prime $t$, a non-negative integer $N$, and any integer $k$, we have
\begin{equation}\label{eq42}
G_{tN+\nu,t}(k\omega^j_t,q) = \dfrac{q^{tk^2+k}}{\displaystyle\prod_{i=0}^{t-1}(\tilde{q};\tilde{q})_{N+k\delta_{i,j-1}-k\delta_{i,j}+\lceil\frac{\nu-i}{t}\rceil}},    
\end{equation}
where $1\le j\le t-1$, $\tilde{q} := q^t$, and $\nu\in\{0,1,2,\ldots,t-1\}$.\\    
\end{theorem}
\par The rest of the proof of Theorem \ref{thm41} is analogous to the proof of Theorem \ref{thm11}.\\

\section{Acknowledgments}
The authors would like to thank Frank Garvan for helping them with his QSERIES and TCORE packages in Maple. The authors would also like to thank George Andrews for his kind interest and Ali K. Uncu for helping them with Figure \ref{fig2}. The authors would also like to thank both the anonymous referees for their helpful comments and suggestions.\\

\section{Declaration of Interest Statement}
The authors declare that they have no known competing financial interests or personal relationships that could’ve appeared to influence the work reported in the submitted paper.\\

\section{Data Availability Statement}
Data sharing is not applicable to this article as no datasets were generated or analysed during the current study.\\

\bibliographystyle{amsplain}


\end{document}